\documentclass[conference]{IEEEtran}

\usepackage{cite}
\usepackage{amsfonts}
\usepackage{amsmath}
\usepackage{graphicx}
\usepackage{booktabs}
\usepackage{mathtools}
\usepackage{enumerate} % item
\usepackage{multirow} % multirow for columns spanning multiple rows
\usepackage{mathrsfs}
\usepackage{ulem}
\usepackage{setspace}
\usepackage{algpseudocode}
\usepackage[ruled, linesnumbered]{algorithm2e}
\usepackage{enumitem}
\usepackage{tabulary}
\usepackage{tabularx}
\usepackage[dvipsnames]{xcolor}
\allowdisplaybreaks[4]
\usepackage{bm}
\usepackage{mathtools}
\usepackage{romannum}
\usepackage{url}
\DeclareMathOperator{\E}{\mathbb{E}}
\usepackage{array}
\newcolumntype{L}{>{\centering\arraybackslash}m{2cm}}

\let\oldnl\nl% Store \nl in \oldnl
\newcommand{\nonl}{\renewcommand{\nl}{\let\nl\oldnl}}% Remove line number for one line]

\makeatletter
\newcommand{\removelatexerror}{\let\@latex@error\@gobble}
\makeatother
\begin{document}
%\title{\textit{\textbf{Resilient Distribution System Restoration Considering Transportable Energy Storage and Soft Open Points}}}
% Adjust the space for float figures
\setlength{\textfloatsep}{3pt} 
\setlength{\abovecaptionskip}{3pt}
\setlength{\belowcaptionskip}{10pt}

% Adjust the space for equations
\setlength{\abovedisplayskip}{-5pt}
\setlength{\belowdisplayskip}{6pt}

\title{Resilient Load Restoration in Microgrids Considering Mobile Energy Storage Fleets: A Deep Reinforcement Learning Approach}
\author{
%	Shuhan~Yao,~\IEEEmembership{Student Member,~IEEE,}
%	Peng~Wang,~\IEEEmembership{Fellow,~IEEE,}
%	Xiaochuan~Liu,~\IEEEmembership{Student Member,~IEEE,}
%	Huajun~Zhang,~\IEEEmembership{Student Member,~IEEE,}
%	and~Tianyang~Zhao,~\IEEEmembership{Member,~IEEE}% <-this % stops a space
	
%	\thanks{Manuscript received November 20, 2017; revised February 27, 2018; accepted April 5, 2018. This work was supported by the Future Resilient System (FRS) at the Singapore-ETH Centre (SEC), which is funded by the National Research Foundation of Singapore (NRF) under its Campus for Research Excellence and Technological Enterprise (CREATE) program.}% <-this % stops a space
%	\thanks{S. H. Yao is with the Interdisciplinary Graduate School, Nan yang Technological University, Singapore 639798 (email: syao002@e.ntu.edu.sg)}% <-this % stops a space
%	\thanks{P. Wang is with the School of Electrical and Electronic Engineering, Nanyang Technological University, Singapore 639798 (e-mail: epwang@ntu.edu.sg).}
%	\thanks{T. Y. Zhao (Corresponding author) is with the Energy Research Institute@NTU, Nanyang Technological University, Singapore 639798 (e-mail: zhaoty@ntu.edu.sg).} 
}

%\author{\IEEEauthorblockN{Shuhan Yao, Huajun Zhang}
%	\IEEEauthorblockA{Institute of Catastrophe Risk Management,\\
%		Interdisciplinary Graduate School,\\
%		Nanyang Technological University\\
%		Singapore\\
%		Email: syao002@e.ntu.edu.sg, \\
%		hzhang031@e.ntu.edu.sg}
%	\and
%	\IEEEauthorblockN{Tianyang Zhao}
%	\IEEEauthorblockA{Energy Research Institute @ NTU,\\
%		Nanyang Technological University\\
%		Singapore\\
%		Email: zhaoty@ntu.edu.sg}
%	\and
%	\IEEEauthorblockN{Xiaochuan Liu, Peng Wang}	
%	\IEEEauthorblockA{School of Electrical and \\ Electronic Engineering,\\
%		Nanyang Technological University\\
%		Singapore\\
%		Email: e160064@e.ntu.edu.sg, \\
%		epwang@ntu.edu.sg
%		}
%
%}
%
%\markboth{IEEE TRANSACTIONS ON SMART GRID}% p

% For conferences
\author{\IEEEauthorblockN{Shuhan Yao, Jiuxiang Gu, Huajun Zhang}
	\IEEEauthorblockA{Interdisciplinary Graduate School,\\
		Nanyang Technological University\\
		Singapore\\}
	
	\and
	\IEEEauthorblockN{Peng Wang, Xiaochuan Liu}
	\IEEEauthorblockA{School of Electrical and \\Electronic Engineering,\\
		Nanyang Technological University\\
		Singapore\\}
	
	\and
 	\IEEEauthorblockN{Tianyang Zhao}
	\IEEEauthorblockA{Energy Research Institute @ NTU,\\
		Nanyang Technological University\\
		Singapore\\}
		
}

\maketitle

\begin{abstract}
%Mobile energy storage system (MESS) is a promising solution to enhance distribution system resilience in terms of mobility and flexibility. This paper proposes a rolling integrated restoration strategy to coordinate the scheduling of MESSs, resource dispatching of microgrids and network reconfiguration of distribution systems. The operation strategy of MESSs is modeled by a multi-layer time-space network technique. The integrated strategy takes into account both line outages in distribution systems and road damages in a transportation network. A rolling optimization framework is adopted to dynamically update subsequent system faults and optimize the restoration strategy accordingly. The optimization problem in each rolling iteration is formulated as a mixed-integer linear programming with various temporal-spatial and operation constraints. The proposed model is verified on an integrated test system with Sioux Falls transportation network and four 33-bus distribution systems. The results indicate that distribution systems with MESSs are more resilient because of the effectiveness of MESS mobility.

Mobile energy storage systems (MESSs) provide mobility and flexibility to enhance distribution system resilience. The paper proposes a Markov decision process (MDP) formulation for an integrated service restoration strategy that coordinates the scheduling of MESSs and resource dispatching of microgrids. The uncertainties in load consumption are taken into account. The deep reinforcement learning (DRL) algorithm is utilized to solve the MDP for optimal scheduling. Specifically, the twin delayed deep deterministic policy gradient (TD3) is applied to train the deep Q-network and policy network, then the well trained policy can be deployed in on-line manner to perform multiple actions simultaneously.
The proposed model is demonstrated on an integrated test system with three microgrids connected by Sioux Falls transportation network. The simulation results indicate that mobile and stationary energy resources can be well coordinated to improve system resilience.

\end{abstract}


\IEEEpeerreviewmaketitle

\section{Introduction}

\IEEEPARstart{R}{ecent} major blackouts caused by extreme events lead to catastrophic consequences for the economy and society \cite{Bie2017}. Load restoration is of paramount importance in resilient smart grids \cite{Wang2016f}. 
%When severe blackouts occur, a variety of local resources, e.g., microgrids and distributed energy resources (energy storage systems, etc.), can be utilized to restore critical loads in distribution systems. Moreover, the emerging mobile energy storage systems (MESSs) \cite{IEEE2019} can provide temporal-spatial mobility and coordinate with stationary local resources for an integrated distribution system restoration.
Great progress has been made in coordinating multiple energy resources to effectively restore electricity supply to critical loads after major blackouts \cite{Xu2018}. Microgrids are well utilized to consolidate stationary energy resources \cite{Chen2017}. Moreover, with the increasing installation of charging/discharging facilities \cite{Yao2018}, microgrids can provide plug-and-play integration of mobile energy storage systems (MESSs) for effective service restoration. 
%MESSs are generally vehicle-mounted container battery energy storage systems with standard interfaces that allow for plug-and-play \cite{IEEE2019}. 
The importance of integrating mobile energy resources into critical load restoration in smart grid has been increasingly recognized in recent studies \cite{Yao2018a, Kim2018, Yao2019}. % Lei2018a
EPRI and Department of Defense of the U.S. initiated a project to demonstrate a containerized grid support storage system featuring utility capacity up to 2 MWh committed to enhancing energy security at military facilities \cite{DepartmentofDefense2016}.
%Reference \cite{Gao2017} formulates a proactive resource allocation scheme of electric buses and transportable batteries to restore critical loads. 
Reference \cite{Che2018} proposes a microgrid-based critical load restoration by adaptively forming microgrids and positioning mobile emergency resources after power disruptions. 
Reference \cite{Lei2018} implements resilient routing and scheduling of mobile power sources via a two-stage framework. 
%A resilient scheme for disaster recovery logistics is proposed in [S Lei, recovery logistics], involving scheduling of repair crews and mobile power source and network reconfiguration. 
%Reference [Yao 2019] presents a joint post-disaster restoration scheme for distribution systems with critical loads by integrating the dynamic scheduling of MESSs. 
%Reference \cite{Yao2019} proposed a rolling horizon-based integrated restoration strategy to effectively restore electricity supply to critical loads by MESSs considering uncertainties. 
However, the optimal scheduling is generally formulated as mixed-integer convex program, which is NP-hard and computationally expensive, in terms of a large number of integer or binary variables in large-scale systems \cite{Lu2019}. In addition, accurate forecast information is necessary in the optimization model \cite{Mocanu2018}.

% These algorithms are difficult to transfer from one scenario to another [Narsin 2019, Definition and evaluation].

%In [Yang liu], a post-event restoration method for power distribution systems with internet data centers s as critical loads is proposed. The emergency operation of data centers and restoration of power distribution systems are jointly optimized to minimize the total utility loss. 
%However, few existing researches have investigated the service restoration in IDCs by integrating the emerging MESS. 

%To leverage the flexibility provided by MESS to facilitate the service restoration in microgrids, more detailed integrated service restoration strategy are needed.

%with consideration of uncertainties in load consumption and damage and repair in both distribution and transportation systems.
%But these studies are proposed for generic critical loads, like hospitals, water stations and other critical infrastructures, without considering particular operation of IDCs. 

%For reliable operation and effective restoration in IDCs, [christopher 2016] presents a optimal investment in IDCs battery storega systems to enhance IDC reliability. Uninterruptible power supply systems are utilized to restore electricity. [Jie Li 2018] studies the optimal operation of IDCs with battery energy storage systems considering islanding scenarios during main grid outages. Workload scheduling and shedding scheme is proposed to restore service by onsite generation resources after power disruptions.

Recent advances in deep reinforcement learning (DRL) give rise to tremendous success in solving challenging decision-making problem \cite{Silver2017, Mnih2015}.
%\cite{Silver2017, Mnih2015, Jaderberg2019}.
 In general, the decision-making problem under uncertainties is formulated using Markov decision process (MDP) \cite{Sheskin2011} and solved iteratively by data-driven DRL algorithms \cite{Mnih2015}.
The application of deep reinforcement learning in energy management systems has been increasingly recognized.
Reference \cite{Foruzan2018} presents a reinforcement learning approach for optimal distributed energy management in a microgrid. A DRL-based economic dispatch in microgrid is proposed in \cite{Liu2018e}. Reference \cite{Xu2019} developed an MDP formulation for the joint bidding and pricing problem and applied DRL algorithm to solve it. Reference \cite{Lu2019} proposes a demand response for home energy management based on DRL. 
%Optimal energy management strategies for energy internet via DRL are demonstrated in \cite{Hua2019}. 
An MDP formulation for electrical vehicle charging is proposed to jointly coordinate a set of charging stations \cite{Sadeghianpourhamami2018}.
A dynamic distribution network reconfiguration using reinforcement learning is proposed in \cite{Gao2019}.
However, research in this area is still in the early stage, the benefit of applying DRL in coordinated scheduling of stationary and mobile energy resources has not yet been fully investigated and further studies are needed. 

%This can deal with uncertainties and be deployed in on-line manner \cite{}.
%
%
%Reinforcement learning can be leveraged to solve the Markov decision process and deal with curse of dimensionality \cite{}.

To address the aforementioned issue, a novel MDP formulation for critical load restoration in microgrids is proposed considering the stationary and mobile energy resources. Uncertainties in load consumption are taken into account. The agent aims to maximize the service restoration in microgrids by jointly coordinating the resource dispatching of microgrids and scheduling of MESS. The MESS fleets are dynamically dispatched among microgrids for load restoration in coordination with microgrid operation.
The proposed model is solved by twin delayed deep deterministic policy gradient (TD3) \cite{Fujimoto2018}, which is an actor-critic algorithm that can deal with discrete or continuous variables in state and action space.
%Reinforcement learning

The remainder of this paper is organized as follows. Section \ref{sec:mathematical modeling} mathematically describes the scheduling of MESSs and integrated service restoration strategy. Section \ref{sec:reinforcement learning} develops the MDP formulation and deep reinforcement learning algorithm. Section \ref{sec:case studies} provides case studies and the paper is concluded in Section \ref{sec:conclusions}.

\section{Mathematical modeling}
\label{sec:mathematical modeling}
%\label{sec: Modeling of Mobile Energy Storage Systems}

%"Markov Decision Process is a framework allowing us to describe a problem of learning from our actions to achieve a goal." \url{https://harderchoices.com/2018/02/11/finite-markov-decision-process-a-high-level-introduction/}

%exogenous information
%
%incrementally gathered data 
%
%finite horizon Markov Decision Process

%\subsection{MDP State Space}
%
%The state space at interval $t$ is represented by a 3-dimensional feature vector 
%$s_t = [t,~p_{\text{d}, t},~\text{soc}_t] \in \mathbf{S}, \forall t \in \mathbf{T}$
%
%
%
%overall state space vector
%
%
%\begin{equation} \label{eq:state_space_1}
%\uline{\text{soc}} \leq \text{soc}_t \leq \overline{\text{soc}}, \forall t \in \mathbf{T}
%\end{equation}
%
%\subsection{MDP Action Space}
%can $a_t$ be a vector? dimensional actions?
%
%Given the state $s_t$, the action $a_t =[p_{\text{g},t}]$ indicates the generation power with constraints as below.
%
%\begin{equation} \label{eq: action_space_1}
%\uline{p}_\text{g} \leq p_{\text{g},t} \leq \overline{p}_\text{g}, \forall t \in \mathbf{T}
%\end{equation}
%
%\subsection{Dynamics}
%The transition will take place as follows.
%
%\begin{equation}
%soc_{t+1} = soc_t - \frac{p_{\text{g},t}}{E_\text{g}}
%\end{equation}
%
%\subsection{MDP Reward}
%
%The reward is as follows.
%
%There are several constraints the system need to comply with.
%
%
%
%\begin{equation}
%p_{\text{g},t} = p
%\end{equation}

%\begin{equation}
%\begin{aligned}
%\min~&-4x_1 - x_2 \notag \\
%\text{s.t.}~&-x_1 + 2x_2 \leq 4, \\
%&2x_1 + 3x_2 \leq 12, \\
%&x_1 - x_2 \leq 3, \\
%&x_1 \geq 0, x_2 \in \textbf{Z}
%\end{aligned}	
%\end{equation}

\subsection{Uncertainties  Modeling}

Uncertainties have been considered including forecasting errors in load consumption. A normal distribution is used to represent the forecasting error of load consumption \cite{Yao2019}. At each time step $t$, the load will be simulated as exogenous information input.

\subsection{Scheduling of Mobile Energy Storage Fleets}

A transportation network is modeled as a weighted graph $\mathcal{G}_\text{T}=(\mathcal{N}_\text{T}, \mathcal{E}_\text{T}, \mathcal{W}_\text{T})$, where $\mathcal{N}_\text{T}$ is the nodes set, 
%representing microgrids and depots' locations. 
while $\mathcal{E}_\text{T}$ denotes the edges set of roads with the edge distance $w  \in \mathcal{W}_\text{T}$.
A set of microgrids $\mathcal{M}$ indexed by $m$ and a set of depots $\mathcal{D}$ are located in the transportation network $\mathcal{G}_\text{T}$.  Location mappings $f_\text{M}:\mathcal{M} \rightarrow \mathcal{N}_\text{T}$ and $f_\text{D}:\mathcal{D} \rightarrow \mathcal{N}_\text{T}$ denote microgrids and depots' locations in the transportation network, respectively. $\Omega$ represents an MESS fleet. An MESS $\omega \in \Omega$ is initially located at a depot $d\in \mathcal{D}$, where it starts and travels among microgrids to provide power supply to power grids, finally it goes back to a depot. 

The scheduling of MESS fleets is defined as a sequence of trips. An MESS $\omega$'s current location at $t$ is represented by $n_\omega^t$, which is generally defined as the node in the transportation network \cite{Yu2019}. In addition, MESS may change destination during its' movement without having to arrive at the next destination, that is,  MESS may be on the edge at $t$, so the location of MESS is defined as $n_\omega^t \in \mathcal{N}_\text{T} \cup \{ (\hat{n}, \delta_{\hat{n}}, \check{n}, \delta_{\check{n}} ) | (\hat{n}, \check{n}) \in E_\text{T}, \delta_{\hat{n}} + \delta_{\check{n}} = w_{\hat{n} \check{n}}, \delta_{\hat{n}} \geq 0, \delta_{\check{n}} \geq 0  \}$, where the $\{ (\hat{n}, \delta_{\hat{n}}, \check{n}, \delta_{\check{n}} ) | (\hat{n}, \check{n}) \in E_\text{T}, \delta_{\hat{n}} + \delta_{\check{n}} = w_{\hat{n} \check{n}}, \delta_{\hat{n}} \geq 0, \delta_{\check{n}} \geq 0  \}$  denotes a location on the edge $(\hat{n}, \check{n}) \in E_\textit{T}$, $\delta_{\hat{n}}$ and $\delta_{\check{n}}$ depict the location's distance to corresponding nodes, and $w_{\hat{n} \check{n}}$ represents the edge length.

The movement decision for MESS $\omega$ at $t$ is to designate the destination $\kappa_\omega^t \in \mathcal{M} \cup \mathcal{D}$, which specifies the destination to one of microgrids or stations. The MESS $\omega$ moves from the current location $n_\omega^t$ and follows the movement decision $\kappa_\omega^t$ to the designated destination. And It is assumed that the MESS $\omega$ always takes the shortest path, which is determined by the Dijkstra's algorithm \cite{Cormen2009}. Therefore, a location function $f_\text{L}$  is defined to obtain the next location $n_\omega^{t+1}$ in graph $G_\text{T}$, by using Dijkstra algorithm based on current location $n_\omega^t$ and designated destination $\kappa_\omega^t$. Thus, we have

\vspace{-6pt}
\begin{gather}
n_\omega^{t+1}  = f_\text{L} (n_\omega^t, \kappa_\omega^t), \forall \omega, t
\end{gather}
\vspace{-16pt}

%$\wp, \kappa$

Binary variables $\zeta_{\omega m}^t$ denote if MESS $\omega$ stays at microgrid $m$ during the interval $t$, which is described as follows.

\vspace{-6pt}
\begin{gather}
\label{eq:mess_zeta}
\zeta_{\omega m}^{t}= 
\begin{cases}
1, & \text{If}~n_\omega^t = n_\omega^{t+1}~\text{and}~n_\omega^t \in f_\text{M} (m) \\
0, & \text{Otherwise} \\
%y & = & \sin(t) \\
\end{cases}
, \forall \omega, m, t
\end{gather}
\vspace{-8pt}

%where the 
%temporal-spatial constraints

%state vector: $s=[t, P_{\text{d},m}^t, Q_{\text{d},m}^t, E_{\text{DG},m}^t, \text{SOC}_\omega^t, \text{CL}_\omega^t, \text{CD}_\omega^t, \text{RT}_\omega^t]$

%state space:

MESS fleets can exchange power with microgrids by charging from or discharging to microgrids. The operation constraints are described as follows.
%operations constraints

\vspace{-6pt}
\begin{gather}
\label{eq:mess_op_1}
\sum_{m \in M} \zeta_{\omega m}^t \leq 1, \forall \omega, t\\
\label{eq:mess_op_2}
-\overline{P}^\omega_\text{ch} \sum_{m \in M} \xi^t_{\omega m} \leq P_\omega^t \leq \overline{P}^\omega_\text{dch} \sum_{m \in M} \xi^t_{\omega m},
\forall \omega, t \\
\label{eq:mess_op_3}
\begin{align}
\text{SOC}_\omega^{t+1}= 
\begin{cases}
\text{SOC}_\omega^t - \frac{\eta_\mathrm{ch}^\omega P_\omega^t}{E_\text{c}^\omega} \Delta_t, & \text{if}~P_\omega^t < 0 \\
\text{SOC}_\omega^t -  \frac{P_\omega^t}{\eta_\text{dch}^\omega E_\text{c}^\omega} \Delta_t, & \text{if}~P_\omega^t \geq 0 \\
%y & = & \sin(t) \\
\end{cases}
, \forall \omega, t
\end{align} \\
\label{eq:mess_op_4}
\uline{\text{SOC}}_\omega \leq \text{SOC}_\omega^t \leq \overline{\text{SOC}}_\omega,
\forall \omega, t
\end{gather}
\vspace{-16pt}

\noindent where $P_\omega^t$ represent the charging/discharging power of MESS $\omega$ from/to microgrid $m$ at interval $t$, negative power depicts that MESS charges from microgrids while positive power means that MESS discharge to microgrid. $\overline{P}^\omega_\text{ch}$ and $\overline{P}^\omega_\text{dch}$ are maximum charging/discharging power of MESS $\omega$. $\text{SOC}_\omega^t$ indicates the state-of-charge (SOC) of MESS $\omega$ at time point $t$. $\uline{\text{SOC}}_\omega$ and $\overline{\text{SOC}}_\omega$ provide the prescribed minimum and maximum level of SOC. $\eta_\text{ch}^\omega$  and $\eta_\text{dch}^\omega$ are charging/discharging efficiency. $E^\omega_c$ indicates the battery capacity of MESS $\omega$.
Constraints \eqref{eq:mess_op_1} indicates that an MESS can only stay at no more than one microgrid, which is also implicated in the Equation \eqref{eq:mess_zeta}. Constraint \eqref{eq:mess_op_2} shows the relation between charging/discharging and temporal-spatial behaviors. That is, only when staying at a microgrid $m$ can MESS $\omega$ charge or discharge to exchange power. Equation \eqref{eq:mess_op_3} calculates the SOC of MESS $\omega$ and Constraint \eqref{eq:mess_op_4} sets the upper and lower bound for SOC. 

%where $0 \rightarrow M-1$ indicates the vehicle stays at corresponding station, while M means the vehicle is on movement. Which represents the remaining time to next destination
%
%action vector: $a = [P_{\text{DG}, m}^t, Q_{\text{DG}, m}^t, P_\omega^t, \text{Dir}_\omega^t]$
%
%action space:
%
%\begin{gather}
%%
%Dir_\omega^t = \{0, 1, \cdots, M-1\} 
%\end{gather}
%

%global state vector
%
%%\begin{figure*}[!tb]
%%\begin{equation}
%%s = [P_{\text{d}, 1}^t \cdots P_{\text{d}, M}^t, Q_{\text{d}, 1}^t \cdots Q_{\text{d}, M}^t, E_{\text{DG},1}^t \cdots E_{\text{DG},M}^t, SOC_1^t, \cdots SOC_\Omega^t, g_\omega^t, k_\omega^t]
%%\end{equation}
%%\end{figure*}
%
%\begin{align}
%s =[& P_{\text{d}, 1}^t \cdots P_{\text{d}, M}^t, Q_{\text{d}, 1}^t \cdots Q_{\text{d}, M}^t, \notag \\
%& E_{\text{DG},1}^t \cdots E_{\text{DG},M}^t, \notag \\
%& SOC_1^t, \cdots SOC_\Omega^t, g_\omega^t, k_\omega^t]
%\end{align}

%\subsubsection{Distribution System Agent}
%~
%
%state action
%
%state space
%
%action vector: $a = [\alpha_{ij}]$ $P_ij \leq $

\subsection{Joint Service Restoration }
\label{subsec:joint_service_restoration}
The operation constraints of microgrids are as follows.

\vspace{-6pt}
\begin{gather} 
\label{eq:joint mg 1}
P_{\text{dg}, m}^t + \sum_{\omega \in \Omega}\xi_{\omega m}^t P_\omega^t = P_{\text{r}, m}^t, 
\forall m, t \\
\label{eq:joint mg 2}
Q_{\text{dg}, m}^t = Q_{\text{r}, m}^t, 
\forall m, t \\
\label{eq:joint mg 3}
0 \leq P_{\text{r}, m}^t \leq P_{\text{load}, m}^t, 
\forall m, t \\
\label{eq:joint mg 4}
Q_{\text{r}, m}^t = P_{\text{r}, m}^t \tan(\cos^{-1} \varphi_m), 
\forall m, t \\
%
%\label{eq: joint_mg_2}
%Q_{\text{dg}, m}^t = Q_{\text{r }, m}^t - Q_{\text{pv}, m}^t, 
%\forall \omega, m, t \\
%
%\label{eq: joint_mg_3}
%0 \leq P_{\text{R}, m}^t \leq P_{\text{D}, m}^t, 
%\forall \omega, m, t \\
%
%\label{eq: joint_mg_4}
%Q_{\text{R}, m}^t = P_{\text{R}, m}^t \tan(\cos^{-1} \varphi_i), 
%\forall \omega, m, t \\
%
\label{eq:joint mg 5}
0 \leq P_{\text{dg}, m}^t \leq \overline{P}_{\text{dg}, m},
\forall m, t \\
\label{eq:joint mg 6}
-\overline{Q}_{\text{dg}, m} \leq Q_{\text{dg}, m}^t \leq \overline{Q}_{\text{dg}, m},
\forall m, t \\
%
%\label{eq: joint_mg_6}
%-\overline{Q}_{\text{DG}, m} \leq Q_{\text{DG}, m}^t \leq \overline{Q}_{\text{DG}, m},
%\forall m, t \\
%
\label{eq:joint mg 7}
E_{\text{dg},m}^{t+1} = E_{\text{dg},m}^{t} - P_{\text{dg},m}^{t+1}\Delta t,
\forall m, t \\
%
%\label{eq: joint_mg_4}
%E_{\text{dg},m}^{t+1} = E_{\text{dg},m}^{t} - P_{\text{dg},m}^{t+1}\Delta t,
%\forall m, t \\
%
\label{eq:joint mg 8}
\uline{E}_{\text{dg}, m}^\text{min} \leq E_{\text{dg},m}^{t} \leq \overline{E}_{\text{dg}, m}
\forall m, t
\end{gather}
\vspace{-16pt}

\noindent where $P_{\text{dg}, m}^t, Q_{\text{dg}, m}^t$ are the active/reactive power generation of equivalent dispatchable DG in microgrid $m$ in interval $t$, respectively. $\overline{P}_{\text{dg}, m}, \overline{Q}_{\text{dg}, m}$ are the maximum active/reactive power generation, respectively.
$P_{\text{r}, m}^t, Q_{\text{r}, m}^t$ are active/reactive load restoration in microgrid $m$, respectively. 
$\varphi_m$ is the power factor.
$E_{\text{dg}, m}^t$ is the energy of equivalent DG. $\overline{E}_{\text{dg}, m}^t$ and $\uline{E}_{\text{dg}, m}^t$ are the energy capacity and minimum energy reserve in microgrid $m$. 
Constraints \eqref{eq:joint mg 1}-\eqref{eq:joint mg 2} describe the active/reactive power balance at microgrid $m$ in interval $t$. It takes into account the power generation of dispatachable DG and mobile energy storage by considering if the location of MESSs. Equations \eqref{eq:joint mg 3}-\eqref{eq:joint mg 4} constrain the load restoration and power factor. Constraints \eqref{eq:joint mg 5}-\eqref{eq:joint mg 6}
depict the power generation capacity. Equation \eqref{eq:joint mg 7} calculates the energy in each microgrid. Constraint \eqref{eq:joint mg 8} presents the upper and lower bounds of energy.

In the wake of major disturbances, the restoration strategy is implemented across multiple microgrids over the horizon to reach a higher level of resilience, which is more focused on the system cost in this work. Therefore, the objective is formulated as follows to minimize the system overall cost.

%\begin{equation}  \label{eq: objective function}
%%\setlength\abovedisplayskip{-5pt}
%\begin{aligned}
%\min~\bigg[\sum_{i\in \mathbf{N}} W_i P_{\mathrm{r},i}^t - \sum_{m\in \mathbf{M}} C_{\mathrm{gen},m}P_{\text{DG},m}^t \\ 
%- \sum_{\omega\in \mathbf{\Omega}} C_{\text{bat},\omega} |P_\omega^t| \bigg] \Delta T
%\end{aligned}
%\end{equation}

\vspace{-6pt}
\begin{align}  
\label{eq:objective function}
\min~\sum_{t \in T} \bigg[ 
\sum_{m \in \mathcal{M}} W_m (P_{\text{load}, m}^t - P_{\text{r}, m}^t) 
+ \sum_{m \in \mathcal{M}} C_{\text{gen},m}P_{\text{dg},m}^t  & \notag \\
+ \sum_{\omega\in \Omega} C_{\text{bat},\omega} |P_\omega^t| + 
\sum_{\omega \in \Omega} C_{\text{tran}, \omega}  (1-\sum_{m \in \mathcal{M}} \zeta_{\omega m}^{t}) 
\bigg] \Delta T &
%V_{\text{avg}, \omega}
%\setlength\abovedisplayskip{-5pt}
%
%\min~\sum_{t \in T} \bigg[ C_\text{iw} \sum_{\psi \in \Psi} (P_{\text{iw}, \psi}^t - \sum_{m \in M} P_{\text{iw}, \psi m }^t ) + C_\text{bw} \sum_{m \in M}  \varTheta_m^{|\mathcal{T}|} & \notag \\
% +  \sum_{m \in \mathbf{M}} C_{\mathrm{gen},m}P_{\text{dg},m}^t  + \sum_{\omega\in \mathbf{\Omega}} C_{\text{bat},\omega} |P_\omega^t| \notag \\ 
% + \sum_{\omega \in \Omega} C_{\text{tran}, \omega} V_{\text{avg}, \omega} \bigg] \Delta T & \notag
\end{align}
\vspace{-14pt}

\noindent where the overall cost is composed of four parts. The first term $\sum_{t \in T}  \sum_{m \in \mathcal{M}} W_m(P_{\text{load}, m}^t - P_{\text{r}, m}^t) $ represents the customer interruption cost. $\sum_{t \in T} \sum_{m \in \mathcal{M}} C_{\text{gen},m}P_{\text{dg},m}^t$ is the microgrids generation cost. The third term $\sum_{t \in T} \sum_{\omega\in \Omega} C_{\text{bat},\omega} |P_\omega^t|$ shows the MESS battery maintenance cost. The last term $\sum_{t \in T} \sum_{\omega \in \Omega} C_{\text{tran}, \omega} V_{\text{avg}, \omega} (1-\sum_{m \in \mathcal{M}} \zeta_{\omega m}^{t})$ calculates the transportation cost of MESSs.

%$C_\text{iw} \sum_{\psi \in \Psi} (P_{\text{iw}, \psi}^t - \sum_{m \in M} P_{\text{iw}, \psi m }^t ) \Delta T $ is the interactive workloads curtailment cost at interval The second term $t$. $C_\text{bw} \sum_{m \in M}  \varTheta_m^{|\mathcal{T}|} \Delta T$ describes the curtailment cost for batch workloads at interval $t$.
%The third term $\sum_{m \in \mathbf{M}} C_{\mathrm{gen},m}P_{\text{dg},m}^t \Delta_t$ shows the microgrid generation cost. 
%The fourth term $\sum_{\omega\in \mathbf{\Omega}} C_{\text{bat},\omega} |P_\omega^t| \Delta T $calculates the battery maintenance cost, while the last term $\sum_{\omega \in \Omega} C_{\text{tran}, \omega} V_{\text{avg}, \omega} \Delta T$ denotes the transportation cost.

\section{Deep Reinforcement Learning Algorithm}
\label{sec:reinforcement learning}

\subsection{Markov Decision Process}

%Q-learning deal with high-dimensional state space, but suffers from curse of dimensionality  in dealing with continuous action space. Thus policy gradient can ****, then to actor-critic methods.
%
%leverage the recent great progress in reinforcement learning  solve Markov decision process.
%MDP offers a standard formalism for describing multi-state? decision making in a probabilistic environment.

The sequential decision-making problem in a stochastic environment is formulated by Markov decision processes (MDPs). In an MDP, an agent observes the state $s_t $ at each time step $t \in \mathcal{T}$ and continually interacts with an environment by following a policy $\pi$ to select actions $a$. In response to the actions, the environment presents new states  $s_{t+1}$ and give rise to rewards $r_t$ to the agent.
An MDP is defined by a 4-tuple $(\mathcal{S}, \mathcal{A}, \mathcal{P}, \mathcal{R})$, where $\mathcal{S}, \mathcal{A}, \mathcal{P}, \mathcal{R}$ are the state space, action space, transition probability functions that satisfy Markov property \cite{Sheskin2011} (i.e., the next state is only dependent on present state and action), and reward functions. 
% Scheduling the operation of a connected vehicular network using deep reinforcement leraning.
The detailed formulation is described as follows.

The state is a vector defined as $s_t = \{ t, P_{\text{load}, m}^t, n_\omega^t, \text{SOC}_\omega^t,E_{\text{dg},m}^{t} \} \in \mathcal{S} $, presenting information on time step, load, the location and SOC of MESSs, and energy in microgrids. 

Furthermore, the action is a vector consisting of decision variables on the designated destination of MESSs charging/discharging behavior of MESSs and generation output in microgrids. The action is defined as $a_t = \{ \kappa_\omega^t, P_\omega^t, P_{\text{dg}, m}^t \} \in \mathcal{A} $. It is noted that $\kappa_\omega^t$ represents categorical action and needs to be one-hot encoded.

%\begin{gather}
%s_t = [pfdas, dfad, dfas] 
%\end{gather}

%\noindent where $dpff = [fdasfd_1, fdsfds_2]$ is ***. 
%
%state variables constraints will be constrained in state transition

%\begin{gather}
%a_t = [ddd_t, dfdas_t]
%\end{gather}

%\noindent where $dd = [ddfd_1..n]$ is ***. Similarly, the categorical actions need to be one-hot encoded.

The state transition $\mathcal{P}: \mathcal{S} \times \mathcal{A} \times \mathcal{S} \rightarrow [0, 1]$ represents the dynamics of the environment,  the transition function indicates the mapping between states at two adjacent time points, so we have $s_{t+1} = \mathcal{P}(s_t, a_t)$. 
To model the uncertainties in load consumption. the exogenous information $P_{\text{load}, m}^t$ in state vector $s_t$ are random variables. Based on the state and action, the next state $s_{t+1}$ can be obtained. In reinforcement learning, the $\mathcal{P}$ is unknown and needs to be learned through interactions between the agent and the environment \cite{Sutton2018}. 
%via Equations \eqref{eq: }-\eqref{eq: } in Section.

%In addition, if the any element in state vector is out of range, as shown in Equations \eqref{eq: }-\eqref{eq: }, the output is corrected by fixing it to either lower or upper bound [OPF, FUZZY]. Suppose, 

%\begin{gather}
%v_t,i = 
%\end{gather} 

%unknown transition function

The reward function is defined as $\mathcal{R}: \mathcal{S} \times \mathcal{A} \times \mathcal{S} \rightarrow \mathbb{R}$, where $r_t = \mathcal{R} (s_t, a_t, s_{t+1})$ is the immediate reward the agent receives by taking action $a_t$ given state $s_t$.
%unknown to the agent.
The immediate reward $r_t$ has two components to take into objectives and penalty violating constraints \cite{Liu2018e}.
The detailed definition is as follows.

%The reward is defined based on the objective and associated constraints 

\vspace{-12pt}
\begin{align}
r_t = \lambda_1 R_{\text{obj}, t} + \lambda_2 C_{\text{pen}, t}
%\begin{cases}
%r_t, \text{If no constraints violated}&  \\
%r_t + \gamma,  ~~~~~~~~~~~~~~~~~~~` \text{Otherwise}&
%\end{cases}
\end{align}
\vspace{-16pt}

\noindent where $\lambda_1, \lambda_2$ are coefficients. $R_{\text{obj}, t} = 
[
\sum_{m \in \mathcal{M}} W_m P_{\text{r}, m}^t
- \sum_{m \in \mathcal{M}} C_{\text{gen},m}P_{\text{dg},m}^t  
- \sum_{\omega\in \Omega} C_{\text{bat},\omega} |P_\omega^t| - 
\sum_{\omega \in \Omega} C_{\text{tran}, \omega} V_{\text{avg}, \omega} (1-\sum_{m \in \mathcal{M}} \zeta_{\omega m}^{t})  ] \Delta T $ 
relates to objective function \eqref{eq:objective function} and is obtained by ignoring the constant term and taking minus sign, thus the cost minimization is transformed into a reward maximization problem. 
The second term $C_{\text{pen}, t}$ is Lagrangian penalty term incurred by violation of constraints. 

%$\lambda_\text{p}$ is the penalty coefficient. The penalty term is related to the constraints \eqref{eq: }-\eqref{eq:}
%
%The rationale behind the definition is that, when a state $sss$ ap, which can prevent the agent from taking actions to violate the system operation constraints
%
%The minus sign is to transfer the cost minimizing optimization problem to a standard reward maximization form in reinforcement learning. 
%
%cumulative penalty term [optimal energy management, Haochen Hua]
%
%The agent will check the magnitudes of all state variables, including local bus voltage, ***. If any of constraints are violated, the reward will be set to ***.
%
%[X. Y. Shang, Discrete reactive power optimization]

\subsection{Twin Delayed Deep Deterministic Policy Gradient}
In reinforcement learning, the return is defined as the sum of discounted reward $G_t = \sum_{i=t}^{T} \gamma^{(i-t)} r(s_i, a_i)$, where $\gamma \in [0, 1]$ is the discount factor.
A policy $\pi: \mathcal{S} \rightarrow \mathcal{A}$ is a mapping from states to selecting actions, i.e., stochastic policy $a_t \sim \pi (\cdot | s_t)$ or deterministic policy $a_t = \pi (s_t)$.
Solving an MDP is to find a policy $\pi$ that maximizes the expected return $\E_{a_t, s_t} [G_0]$.

In order to deal with continuous and discrete variables in state and action space, an actor-critic algorithm is adopted \cite{Lillicrap2015}, e.g. deep deterministic policy gradient (DDPG) and twin delayed deep deterministic policy gradient (TD3), which concurrently learn Q-functions and a policy. It uses off-policy data and the Bellman equation to learn the Q-function, and uses the Q-function to learn the policy $\pi_{\phi}$ parameterized with $\phi$\cite{Sutton2018}.

%The environment is assumed to be fully observed by agent.
%
%An agent 
%A deterministic actor function is $a=\pi_{\phi}(s)$ parameterized by $\phi$.

%In Q-learning, the action-value function $Q_\pi (s_t, a_t) = \E_\pi [G_t | s_t, a_t] $ denotes the expected return starting from $s_t$, taking the action $a_t$ and thereafter following policy $\pi$, 

%the following action-value function
%\begin{align}
%Q_\pi (s_t, a_t) = \E_\pi [G_t | s_t, a_t] = \E_\pi [\sum_{i=t}^{T} \gamma^{(i-t)} r(s_i, a_i) ]
%\end{align}
%
%\noindent where $\E_\pi[\cdot] $ depicts the expected values of random variables given the policy $\pi$, and the $Q_\pi$ is the action-value function for policy $\pi$
%The optimal Q function obeys the Bellman equation, which describes the optimal $Q^*(s_t, a_t)$ in recursive form, $Q^* (s_t, a_t) = \E_{s_{t+1}, a_{t+1}}[ r_t +  \gamma  \max_{a_{t+1}} Q^* (s_{t+1}, a_{t+1}) ]$ \cite{Mnih2015}.

%\begin{align}
%Q^* (s_t, a_t) = \E_{s_{t+1}, a_{t+1}}[ r_t +  \gamma  \max_{a_{t+1}} Q^* (s_{t+1}, a_{t+1}) ]
%\end{align}

In Q-learning of TD3, the Q-function is estimated by a differentiable function approximator $Q_\theta (s, a)$, which is a neural network with weights $\theta$ as a Q-network.  
The Q-network can be learned to reduce the mean-squared Bellman error. To make the training converge and stable, a separate target Q-network $Q_{\theta'}$ and a target policy network $\pi_{\phi'}$ are utilized to generate optimal target value \cite{Mnih2015}.
%where the $r_t + \gamma \max_{a_{t+1}} Q^* (s_{t+1}, a_{t+1} ) $ is substituted with approximate target values $r_t + \gamma \max_{a_{t+1}} Q_\theta' (s_{t+1}, \pi{s_{}} a_{t+1})$ \cite{Lillicrap2015}. 
Therefore a sequence of loss functions $L(\theta)$ is set up by the mean-squared Bellman error as $L(\theta) = \E_{s_t, a_t, r_t} [(Q_\theta (s_t, a_t) - y_t)^2]$, where the target value $y_t$ is defined as $y_t = r_t + \gamma \max_{a_{t+1}} Q_{\theta'} (s_{t+1}, \pi_{\phi'}(s_{t+1}))$.

%\begin{align}
%\label{eq: learning_loss}
%L(\theta) = \E_{s_t, a_t, r_t} [(Q_\theta (s_t, a_t) - y_t)^2]
%\end{align}

%\noindent where the target value $y_t$ is defined as:
%
%\begin{gather}
%y_t = r_t + \gamma \max_{a_{t+1}} Q_{\theta'} (s_{t+1}, \pi_{\phi'}(s_{t+1}))
%\end{gather} 

%\begin{align}
%\label{eq: learning_target}
%y_t = \left\{
%\begin{aligned}
% r_t, \hfill \text{if episode terminates at step}~t & \\
% r_t + \gamma \max_{a_{t+1}} Q_\theta (s_{t+1}, a_{t+1} ), ~~~~~\hfill \text{otherwise} & \\
%%y & = & \sin(t) \\
%\end{aligned}
%\right.
%\end{align}

%\begin{align}
%y_t = 
%\begin{cases}
% r_t,  \hfill \text{If episode terminates at step}~t &  \\
% r_t + \gamma \max_{a_{t+1}} Q_\theta (s_{t+1}, a_{t+1} ),  ~~~\text{Otherwise}  &
%\end{cases}
%\end{align}

\noindent 
%when $s_{t+1}$ is a terminal state, the target value indicates that the agent gets no additional rewards after current state. [nature 2015]. 

The Q-network is updated by one step gradient descent  using $L(\theta)$. A soft target update is used for actor-critic algorithm \cite{Lillicrap2015}, the target networks are updated by Polyak averaging, $\theta' = \tau \theta + (1-\tau) \theta', \phi' = \tau \phi + (1-\tau) \phi'$, where $\tau \in [0, 1]$ is the Polyak hyperparameter (usually $\tau \ll 1$).
Furthermore, TD3 concurrently learns two Q-networks, $Q_{\theta_1}$ and $Q_{\theta_2}$ by minimizing mean-squared Bellman error. By upper-bounding the less biased value approximator $Q_{\theta_2}$ with the biased estimate $Q_{\theta_1}$, a single target update for clipped Double Q-learning is  
obtained by taking the minimum between the two Q-networks:

\vspace{-8pt}
\begin{gather}
\label{eq: td3_1}
y_t = r_t + \gamma \min_{i=1,2} Q_{\theta'_i} (s_{t+1}, \pi'_{\phi'}(s_{t+1}) )
\end{gather}
\vspace{-13pt}

\noindent Then $Q_{\theta_1}$ and $Q_{\theta_2}$ are updated by minimizing the corresponding mean-squared Bellman error as follows.

\vspace{-6pt}
\begin{gather}
\label{eq:bellman error 2}
L(\theta_i) = \E_{s_t, a_t, r_t} [(Q_{\theta_i} (s_t, a_t) - y_t)^2], \forall i=1,2
\end{gather}
\vspace{-16pt}

Target smoothing regularization is to add a small amount of random noise $\tilde{\epsilon}$ to the target policy network in target update and averaging over mini-batches. The modified target actions $\tilde{a}$ and target values $y_t$ are as follows.

\vspace{-8pt}
\begin{gather}
\label{eq:target action}
\tilde{a}_{t+1} = \pi_{\phi'} (s_{t+1}) + \tilde{\epsilon},~\tilde{\epsilon} \sim \text{clip} (\mathcal{N}(0, \tilde{\sigma}^2), -c, c) \\
\label{eq:target value}
y_t = r_t + \gamma \min_{i=1,2} Q_{\theta'_i} (s_{t+1}, \tilde{a}_{t+1})
\end{gather}
\vspace{-13pt}

%When updating the critics $Q_{\theta_1}$ and $Q_{\theta_2}$, a learning target using a deterministic  policy is highly susceptible to inaccuracies induced by function approximation error, increasing the variance of the target value. This induced variance can be  reduced through regularization by target policy smoothing. Since deterministic polices can overfit to narrow peas in the value estimate, the Q-function over similar actions are smoothed out to have similar value, by adding a small amount of random noise $\tilde{\epsilon}$ to the target policy network in target update and averaging over minibatches. The modified target actions $\tilde{a}$ and target values $y_t$ are thus:
%
%\begin{gather}
%\label{eq:target action}
%\tilde{a}_{t+1} = \pi'_{\phi'} (s_{t+1}) + \tilde{\epsilon},~\tilde{\epsilon} \sim \text{clip} (\mathcal{N}(0, \tilde{\sigma}^2), -c, c) \\
%\label{eq:target value}
%y_t = r_t + \gamma \min_{i=1,2} Q'_{\theta'_i} (s_{t+1}, \tilde{a})
%\end{gather}

\noindent where the added noise is a normal distribution with zero-mean and standard deviation $\tilde{\sigma}$, and clipped by a hyperparameter $c$.

%\vspace{6pt}
%%%%%%%%%%%%%%%%%%%%
%\noindent \textit{3) Delayed Policy Updates}
%\vspace{6pt}
The Policy learning of TD3 is to find a policy $\pi_\phi (s_t)$ that maximizes the expected discounted return $J(\phi) = \E_{s_t, a_t} [G_0] \approx \E_{s_t} [Q_\theta (s_t, \pi_\phi(s_t))]$ \cite{Fujimoto2018}. The policy network is updated by applying the chain rule to the $J(\phi)$ with respect to the actor parameters $\phi$ and gradient ascent is implemented. 
%Thus, the policy learning algorithm is $\nabla_\phi J(\phi) = \E_{s_t} [\nabla_a Q_\theta (s_t, a_t) |_{a_t = \pi_{\phi} (s_t)} \nabla_\phi \pi_{\phi} (s_t) ]$.
The policy $\pi_\phi$ is optimized with respect to $Q_{\theta_1}$ to maximize the expected return $J(\phi)$, so the policy learning is written as:

\vspace{-9pt}
\begin{gather}
\label{eq:policy learning 2}
\nabla_\phi J(\phi) = \E_{s_t} [\nabla_a Q_{\theta_1} (s_t, a_t) |_{a_t = \pi_{\phi} (s_t)} \nabla_\phi \pi_{\phi} (s_t) ] 
\end{gather}
\vspace{-13pt}

In addition, the policy network $\pi_\phi$ is updated at a lower frequency than the value network $Q_{\theta_1}$, in order to reduce error before introducing a policy update \cite{Fujimoto2018}. 

\section{Case Studies}
\label{sec:case studies}

The case studies are implemented on an integrated test system, based on Sioux Falls transportation network 
%\cite{LeBlanc1975} 
and three microgrids, to verify the effectiveness of the proposed service restoration strategy.

\subsection{Test Systems}

\begin{figure}[!tbp]
	\centering
	\includegraphics[width=0.82\columnwidth, clip]{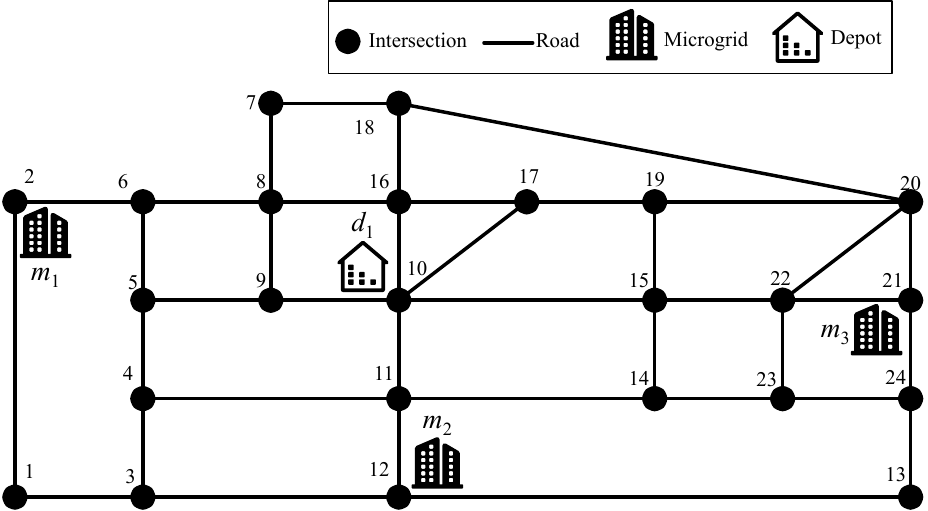}
	\caption{An integrated test system with a Sioux Falls transportation network connecting three microgrids.}
	\label{fig:integrated_test_system}
\end{figure}

Fig. \ref{fig:integrated_test_system} shows an integrated test system with microgrids connected by the Sioux Falls transportation network. 
%The distance of each edge in transportation network is double.
The length of the entire time horizon $T_\text{H}$ is set to 24-h and the length of interval is 1-h. 
A depot is located at node \#10 in the transportation network.
There are three microgrids located at nodes \#2, \#12, \#21 in the transportation network, respectively.  
The operational parameters of microgrids are shown in Table \ref{table:microgrids parameters}.
The predicted value of industrial, commercial and residential loads, as well as prediction intervals could be obtained in \cite{Yao2019}. The parameters for MESS refers to \cite{Yao2019}. The customer interruption cost for industrial, commercial and residential loads are $\$8/\text{kWh}$, $\$10/\text{kWh}$ and $\$2/\text{kWh}$, respectively.
The unit generation cost in microgrid is $\$0.5/\text{kWh}$. The unit battery maintenance cost is $\$0.2/\text{kWh}$. The unit transportation cost is $\$80/\text{h}$.

\begin{table}[!btp]
	\centering
	\caption{Generation Resources and Local Loads for microgrids} \label{table:microgrids parameters}
	\resizebox{0.73\columnwidth}{!}{
		\begin{tabular}{ccccc}
			\hline
			\multicolumn{2}{c}{Microgrid \#}             & 1 & 2 & 3  \\ \hline
			\multirow{4}{*}{Generation} & $\overline{P}_{\text{dg},m}$ (MW) & 1.0 & 1.80 & 1.20  \\ \cline{2-5}
			& $\overline{Q}_{\text{dg}, m}$ (MVar)     & 0.8 & 1.5 & 1.0   \\ \cline{2-5}
			& $\overline{E}_{\text{dg}, m}$ (MWh)     & 20 & 35 & 23   \\ \cline{2-5}
			& $\uline{E}_{\text{dg}, m}$ (MWh)  &   2.0 & 3.5 & 2.3    \\ \hline
			
			\multirow{3}{*}{Load} & Peak load (MW)   &   3.0 & 3.0 & 3.0  \\ \cline{2-5}
			& Power factor &  0.9 & 0.9  & 0.9    \\ \cline{2-5}
			& Load type   &  C & R & I   \\  \hline
			
%			\multicolumn{6}{l}{Note: C - commercial, R - residential, I - industrial}
		\end{tabular}
	}
\end{table}

%% the idle power consumption of a server can be as low as 50-65% of the peak power consumption, which can range from 100-250W [Qureshi 2009-5, 7, 8]. 
%
%Workload traces are adopted in
%
%Predicted interactive workloads , GWA-T-12 trace [Siqi Shenan -2015] and batch workloads are modified according to [jie li-2018], the workload shedding cost is $\$1000/\text{request}$  for both interactive workload and batch workloads.
%%and $\$1000/\text{request}$ for batch workloads.

\subsection{Simulation Results}

 The total cost is $\$351172$, with the customer interruption cost $\$309859$, microgrid generation cost $\$34872$, MESS generations cost $\$4600$ and transportation cost $\$1840$. The load restorations in three microgrids are $72.92\%, 56.66\%$ and $50.76\%$, respectively.

\begin{figure}[!tbp]
	\centering
	\includegraphics[width=\columnwidth, clip]{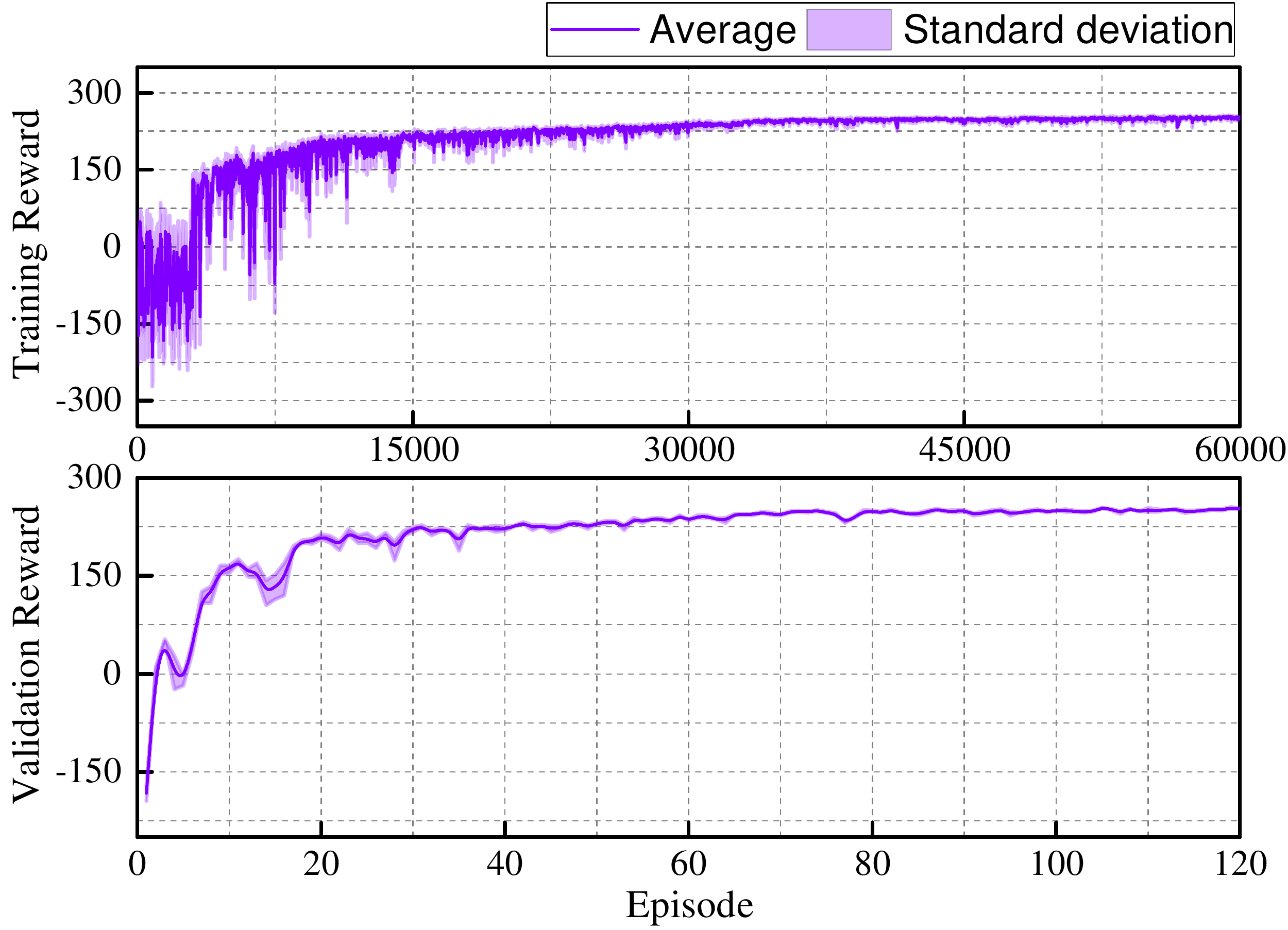}
	\caption{Learning and validation curves.}
	\label{fig:learning curve}
\end{figure}

Fig. \ref{fig:learning curve} illustrates the evolution of learning and validation rewards over 60000 episodes. 
%In order to remove the dependency on the initial parameters of the policy, 
A purely exploratory policy is carried out for the first 3000 episodes. Then, an off-policy exploration strategy is adopted with Gaussian noise.
%with a temporally-correlated noise Ornstein-Hulenbeck process \cite{Uhlenbeck1930}.
In the learning curve, the average and standard deviation are obtained every 10 episodes. In the validation curve, the validation is  evaluated every 500 episodes over 20 episodes with no exploration noise.
It can be seen that the learning process converges to a suboptimal policy in 40000 episodes.
The results indicate that the proposed approach can learn a policy to maximize the cumulative rewards. After learning, the model can be deployed in on-line manner.
%Learning curve

\begin{figure}[!tbp]
	\centering
	\includegraphics[width=\columnwidth, clip]{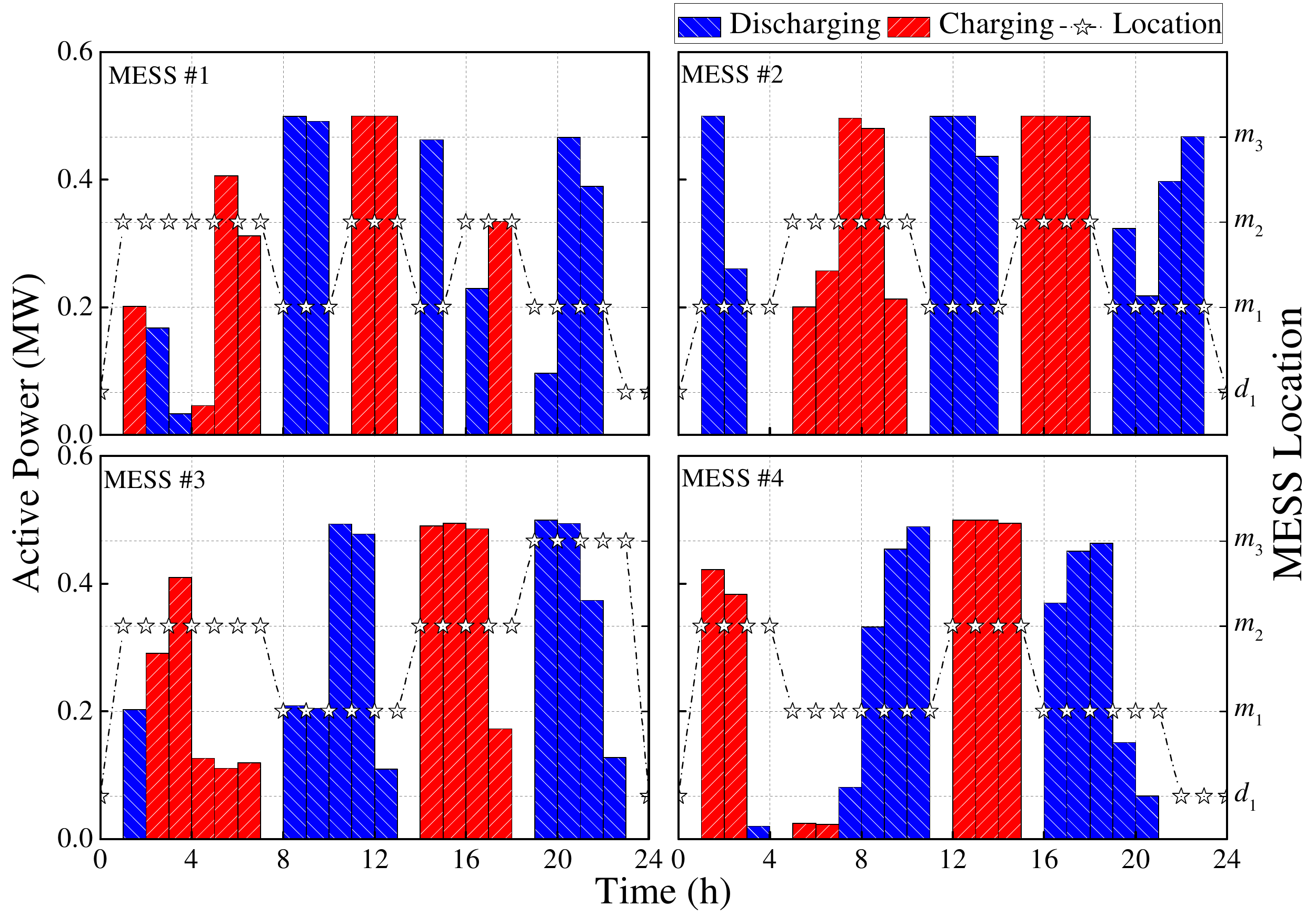}
	\caption{Scheduling results of the MESS fleets.}
	\label{fig:mess scheduling}
\end{figure}

Fig. \ref{fig:mess scheduling} presents the charging/discharging schedule with respect to the position of MESS. The bar shows the charging/discharging active power while the dash lines with asterisks and right Y-axis indicates the MESS's movements. The dynamic scheduling of MESS optimizes the trip chain of MESSs and corresponding charging/discharging behaviors.

The simulation result shows that MESSs transport energy among microgrids to restore critical loads by charging from some microgrids and discharging to others. For example, it is observed that MESS \#1 is dispatched between microgrid \#1 and microgrid \#2. The MESS \#1 initially moves to microgrid \#2 from depot and charges at microgrid \#2. Next, it moves back and forth between microgrid \#2 and microgrid \#1 in (07:00-22:00) to transfer energy. The integration of MESSs and coordination with microgrids can leverage the MESSs mobility. Also, the MESSs can carry out load shifting within the same microgrid. For instance, MESS \#3 charges at microgrid \#2 in (01:00-02:00) and discharges in (02:00-07:00). The results highlight the importance of effective utilization of MESSs mobility and flexibility.

\section{Conclusions}
\label{sec:conclusions}
This paper presents a novel MDP formulation for service restoration strategy in microgrids by coordinating the scheduling of MESSs and resource dispatching of microgrids. The DRL algorithms are leveraged to solve the formulated sequential decision-making problem with consideration of uncertainties in load consumption. The well trained policy can be deployed in on-line manner and is computationally efficient. The simulation results verify the effectiveness of MESSs mobility that transport energy among microgrids to facilitate load restoration. Mobile and stationary resources can be jointly coordinated to enhance system resilience.

\section*{Acknowledgments}
This work was supported by the Future Resilient Systems (FRS) at the Singapore-ETH Centre (SEC).
%,which is funded by the National Research Foundation of Singapore (NRF) under its Campus for Research Excellence and Technological Enterprise (CREATE) program.

\ifCLASSOPTIONcaptionsoff
  \newpage
\fi



%\bibliographystyle{IEEEtran}
% To delete underline in publications' titles.
\normalem
\bibliographystyle{IEEEtran}
\bibliography{20pes_ref.bib}
%bibliography{"/Users/yaoshuhan/OneDrive Business/OneDrive - Nanyang Technological University/Research Program/Mendeley library/BibTeX/lib_ref.bib"}
%\printbibliography

%\newpage


%\vfill


\end{document}